\documentclass[12pt]{amsart}
\usepackage{amsfonts}

\usepackage{amsmath}
\usepackage{amsthm}
\usepackage{graphicx}

\theoremstyle{plain}
\newtheorem{thm}{Theorem}[section]
\newtheorem{prop}[thm]{Proposition}
\newtheorem{cor}[thm]{Corollary}
\newtheorem{lem}[thm]{Lemma}
\theoremstyle{remark}
\newtheorem{rem}{Remark}
\theoremstyle{definition}
\newtheorem{ex}[thm]{Example}

\numberwithin{equation}{section}
\numberwithin{figure}{section}
\numberwithin{table}{section}

\email{inoue@math.sci.hokudai.ac.jp}
\email{nakano\_y@math.sci.hokudai.ac.jp}
\email{v.anh@fsc.qut.edu.au}
\thanks{This work is partially supported by the Australian 
Research Council grant DP0345577.}

\begin{document}
\title[\null]{Linear filtering of systems with memory}
\author[\null]{Akihiko Inoue, Yumiharu Nakano and Vo Van Anh}
\address{Department of Mathematics \\
Faculty of Science \\
Hokkaido University \\
Sapporo 060-0810 \\
Japan}
\address{Department of Mathematics \\
Faculty of Science \\
Hokkaido University \\
Sapporo 060-0810 \\
Japan}
\address{School of Mathematical Sciences \\
Queensland University of Technology \\
GPO Box 2434, Brisbane, Queensland 4001 \\
Australia}

\begin{abstract}
We study the linear filtering problem for systems driven by continuous
Gaussian processes $V_1$ and $V_2$ with memory described by two parameters.
The processes $V_j$ have the virtue that they possess stationary increments
and simple semimartingale representations simultaneously. It allows for
straightforward parameter estimations. After giving the semimartingale
representations of $V_j$ by innovation theory, we derive Kalman-Bucy-type
filtering equations for the systems. We apply the result to the optimal
portfolio problem for an investor with partial observations. We illustrate
the tractability of the filtering algorithm by numerical implementations.
\end{abstract}

\maketitle

\section{Introduction}

\label{sec:1}

In this paper, we use the following Gaussian process $(V(t))_{t\in 
\mathbf{R}%
}$ with stationary increments as the driving noise process: 
\begin{equation}
V(t)=W(t)-\int_{0}^{t}\left( \int_{-\infty 
}^{s}pe^{-(q+p)(s-u)}dW(u)\right)
ds\qquad (t\in \mathbf{R}),  \label{eq:1.1}
\end{equation}
where $p$ and $q$ are real constants such that 
\begin{equation}
0<q<\infty ,\quad -q<p<\infty ,  \label{eq:1.2}
\end{equation}
and $(W(t))_{t\in \mathbf{R}}$ is a one-dimensional Brownian motion
satisfying $W(0)=0$. The parameters $p$ and $q$ describe the memory of $%
V(\cdot )$. In the simplest case $p=0$, $V(\cdot )$ is reduced to the
Brownian motion, i.e., $V(\cdot )=W(\cdot )$.

It should be noticed that (\ref{eq:1.1}) is not a semimartingale
representation of $(V(t))_{0\leq t\leq T}$ with respect to the natural
filtration $\mathcal{F}^{V}(\cdot )$ of $(V(t))_{0\leq t\leq T}$ since $%
(W(t))_{0\leq t\leq T}$ is not $\mathcal{F}^{V}(\cdot )$-adapted. Using the
innovation theory as described in Liptser and Shiryayev \cite{LS} and a
result in Anh et al.~\cite{AIK}, we show (Theorem \ref{thm:2.1}) that there
exists a one-dimensional Brownian motion $(B(t))_{0\leq t\leq T}$, called
the \textit{innovation process}, satisfying 
\begin{equation*}
\sigma (B(s):0\leq s\leq t)=\sigma (V(s):0\leq s\leq t)\qquad (0\leq t\leq
T),
\end{equation*}
and 
\begin{equation}
V(t)=B(t)-\int_{0}^{t}\left( \int_{0}^{s}l(s,u)dB(u)\right) ds\qquad (0\leq
t\leq T),  \label{eq:1.3}
\end{equation}
where $l(t,s)$ is a Volterra kernel given by 
\begin{equation}
l(t,s)=pe^{-(p+q)(t-s)}\left\{ 
1-\frac{2pq}{(2q+p)^{2}e^{2qs}-p^{2}}\right\}
\quad (0\leq s\leq t\leq T).  \label{eq:1.4}
\end{equation}
With respect to the natural filtration $\mathcal{F}^{B}(\cdot )$ of 
$B(\cdot
)$, which is equal to $\mathcal{F}^{V}(\cdot )$, (\ref{eq:1.3}) gives the
semimartingale representation of $V(\cdot )$. Thus the process $V(\cdot )$
has the virtue that it possesses the property of a stationary increment
process with memory and a simple semimartingale representation
simultaneously. We know no other process with this kind of properties. The
two properties of $V(\cdot )$ become a great advantage, for example, in its
parameter estimation.

In \cite{AI} and Anh et al.~\cite{AIK, AIN}, the process $V(\cdot )$ is 
used
as the driving process for a financial market model with memory. The
empirical study for S\&P 500 data in Anh et al.~\cite{AIP} shows that the
model captures very well the memory effect when the market is stable. The
work in these references suggests that the process $V(\cdot )$ can serve as
an alternative to Brownian motion when the random disturbance exhibits
dependence between different observations.

In this paper, we are concerned with the filtering problem of the
two-dimensional partially observable process $(X(t),Y(t))_{0\le t\le T}$
governed by the following linear system of equations: 
\begin{gather}  \label{eq:1.5}
\left\{ 
\begin{split}
dX(t)&= \theta X(t)dt + \sigma dV_1(t), \quad X(0)= X_0, \\
dY(t)&= \mu X(t)dt + dV_2(t), \quad Y(0)=0.
\end{split}
\right.
\end{gather}
Here $X(\cdot)$ and $Y(\cdot)$ represent the state and the observation
respectively. For $j=1,2$, the noise process $V_j(\cdot)$ is described by 
(%
\ref{eq:1.1}) with $(p,q)=(p_j,q_j)$ and $W(\cdot)=W_j(\cdot)$. We assume
that the Brownian motions $W_1(\cdot)$ and $W_2(\cdot)$, whence 
$V_1(\cdot)$
and $V_2(\cdot)$, are independent. The coefficients $\theta, \sigma, 
\mu\in%
\mathbf{R}$ with $\mu\neq 0$ are known constants, and $X_0$ is a centered
Gaussian random variable independent of $(V_1,V_2)$.

The basic filtering problem for the linear model (\ref{eq:1.5}) with memory
is to calculate the conditional expectation $\mathbb{E}[X(t)|\mathcal{F}%
^{Y}(t)]$, called the (optimal) \textit{filter}, where $\mathcal{F}%
^{Y}(\cdot )$ is the natural filtration of the observation process $Y(\cdot 
)
$. In the classical Kalman-Bucy theory (see Kalman \cite{K}, Kalman and 
Bucy 
\cite{KB}, Bucy and Joseph \cite{BJ}, Davis \cite{D} and \cite{LS}),
Brownian motion is used as the driving noise. Attempts have been made to
resolve the filtering problem of dynamical systems with memory by replacing
Brownian motion by other processes. In Kleptsyna et al.~\cite{KKA, KLR, KL}
and others, fractional Brownian motion was used. Notice that fractional
Brownian motion is not a semimartingale. In the discrete-time setting,
autoregressive processes are used as driving noise (see, e.g., \cite{K}, 
\cite{BJ} and Jazwinski \cite{J}). Our model may be regarded as a
continuous-time analogue of the latter since it is shown in \cite{AI} that 
$%
V(\cdot )$ is governed by a continuous-time AR$(\infty )$-type equation.

The Kalman-Bucy filter is a computationally tractable scheme for the 
optimal
filter of a Markovian system. We aim to derive a similar effective 
filtering
algorithm for the system (\ref{eq:1.5}) which has memory. However, rather
than considering (\ref{eq:1.5}) itself, we start with a general continuous
Gaussian process $(X(t))_{0\leq t\leq T}$ as the state process and $Y$
defined by 
\begin{equation*}
Y(t)=\int_{0}^{t}\mu (s)X(s)ds+V(t)
\end{equation*}
as the observation process, where $\mu (\cdot )$ is a deterministic 
function
and $V(\cdot )$ is a process which is independent of $X(\cdot )$ and given
by (\ref{eq:1.1}). Using (\ref{eq:1.3}) and (\ref{eq:1.4}), we derive
explicit Volterra integral equations for the optimal filter (Theorem \ref
{thm:3.1}). In the special case (\ref{eq:1.5}), the integral equations are
reduced to differential equations, which gives an extension to Kalman-Bucy
filtering equations (Theorem \ref{thm:3.4}). Due to the non-Markovianness 
of
the formulation (\ref{eq:1.5}), it is expected that the resulting filtering
equations would appear in the integral equation form (cf. Kleptsyna et al. 
\cite{KKA}). The fact that good Kalman-Bucy-type differential equations can
be obtained here is due the special properties of (\ref{eq:1.5}). This
interesting result does not seem to hold for any other formulation where
memory is inherent.

We apply the results to an optimal portfolio problem in a partially
observable financial market model. More precisely, we consider a stock 
price
model that is driven by the solution $V(\cdot )$ to (\ref{eq:1.1}). 
Assuming
that the investor can observe the stock price but not the drift process, we
discuss the portfolio optimization problem of maximizing the expected
logarithmic utility from terminal wealth. To solve this problem, we make 
use
of our results on filtering to reduce the problem to the case where the
drift process is adapted to the observation process. We then use the
martingale methods\/ (cf.~Karatzas and Shreve \cite{KS}) to work out the
explicit formula for the optimal portfolio (Theorem \ref{thm:4.1}).

This paper is organized as follows. In Section \ref{sec:2}, we prove the
semimartingale representation (\ref{eq:1.3}) with (\ref{eq:1.4}) for $%
V(\cdot )$. Section \ref{sec:3} is the main body of this paper. We derive
closed form equations for the optimal filter. In Section \ref{sec:4}, we
apply the results to finance. In Section \ref{sec:5}, we illustrate the
advantage of $V(\cdot )$ in parameter estimation. Some numerical results on
Monte Carlo simulation are presented. Finally, in Section \ref{sec:6}, we
numerically compare the performance of our filter with the Kalman-Bucy
filter in the presence of memory effect.

\section{Driving noise process with memory}

\label{sec:2}

Let $T\in (0,\infty)$, and let $(\Omega,\mathcal{F},\mathbb{P})$ be a
complete probability space. For a process $(A(t))_{0\le t\le T}$, we denote
by $\mathcal{F}^A(t)$ the $\mathbb{P}$-augmentation of the filtration $%
\sigma(A(s): 0\le s\le t)$, $0\le t\le T$.

Let $V(\cdot)$ be the process given by (\ref{eq:1.1}). The process $%
(V(t))_{t\in\mathbf{R}}$ is a continuous Gaussian process with stationary
increments. The aim of this section is to prove (\ref{eq:1.3}) with (\ref
{eq:1.4}).

By \cite[Theorem 7.16]{LS}, there exists a one-dimensional Brownian motion 
$%
(B(t))_{0\le t\le T}$, called the innovation process, such that 
$\mathcal{F}%
^B(t)=\mathcal{F}^V(t)$, $0\le t\le T$ and that 
\begin{gather}  \label{eq:2.1}
V(t)=B(t)-\int_{0}^{t}\alpha(s)ds, \qquad (0\le t\le T), \\
\alpha(t)=\mathbb{E}\left[\left.\int_{-\infty}^t 
pe^{-(q+p)(t-u)}dW(u)\right|%
\mathcal{F}^V(t)\right] \qquad (0\le t\le T).  \notag
\end{gather}
Thus, $(V(t))_{0\le t\le T}$ is a $\mathcal{F}^V$-semimartingale with
representation (\ref{eq:2.1}).

To find a good representation of $\alpha(\cdot)$, we recall the following
result from \cite[Example 5.3]{AIK}: 
\begin{equation*}
\alpha(t)=\int_{0}^{t}k(t,s)dV(s) \qquad(0\le t\le T)
\end{equation*}
with 
\begin{equation*}
k(t,s)=p(2q+p)\frac{(2q+p)e^{qs}-pe^{-qs}}{(2q+p)^2e^{qt}-p^2e^{-qt}}
\qquad(0\le s\le t).
\end{equation*}
From the theory of Volterra integral equations, there exists a function $%
l(t,s)\in L^2[0,T]^2$, called the resolvent of $k(t,s)$, such that, for
almost every $0\le s\le t\le T$, 
\begin{gather}
\begin{split}
&l(t,s)-k(t,s)+\int_{s}^{t}l(t,u)k(u,s)du=0, \\
&l(t,s)-k(t,s)+\int_{s}^{t}k(t,u)l(u,s)du=0
\end{split}
\label{eq:2.2}
\end{gather}
(see \cite[Chapter 4, Section 3]{D} and \cite[Chapter 9]{GLS}). Using 
$l(t,s)
$, we have the following representation of $\alpha$ in terms of the
innovation process $B$: 
\begin{equation}  \label{eq:2.3}
\alpha(t)= \int_{0}^{t}l(t,s)dB(s) \qquad (0\le t\le T).
\end{equation}
We shall solve (\ref{eq:2.2}) explicitly to obtain $l(t,s)$.

\begin{thm}
\label{thm:2.1} The expression $(\ref{eq:1.4})$ holds.
\end{thm}

\begin{proof}
We have $k(t,s)=a(t)b(s)$ for $0\le s\le t$, 
where, for $t\in [0,T]$,
\[
a(t)=\frac{p(2q+p)}{(2q+p)^2e^{qt}-p^2e^{-qt}},\quad 
b(t)=(2q+p)e^{qt}-pe^{-qt}.
\]
Fix $s\in [0,T]$ and define $x(t)=x_s(t)$ and 
$\lambda=\lambda_s$ by
\[
x(t)=\int_{s}^{t}b(u)l(u,s)du\quad(s\le t\le T),\quad \lambda=b(s).
\]
Then, from (\ref{eq:2.2}) we obtain
\[
\frac{dx}{dt}(t)+a(t)b(t)x(t)=\lambda a(t)b(t),\quad x(s)=0.
\]
The solution $x$ is given by
\[
x(t)=\lambda-\lambda e^{-\int_{s}^{t}a(u)b(u)du},
\]
whence $l(t,s)$ is obtained as
\[
l(t,s)=a(t)b(s)e^{-\int_{s}^{t}a(u)b(u)du}
=k(t,s)e^{-\int_{s}^{t}k(u,u)du}\quad(0\le s\le t).
\]

We have
\[
k(u,u)=p-\frac{2p^2q}{(2q+p)^2e^{2qu}-p^2}.
\]
By the change of variable $x(u)=(2q+p)^2e^{2qu}-p^2$, we obtain
\begin{equation*}
\begin{split}
&2p^2q\int_{s}^{t}\frac{du}{(2q+p)^2e^{2qu}-p^2}
=2p^2q\int_{x(s)}^{x(t)}\frac{1}{2qx(x+p^2)}dx\\
&\qquad=2p^2q\int_{x(s)}^{x(t)}\left\{\frac{1}{x}-\frac{1}{x+p^2}\right\}dx
=\log \frac{x(t)(x(s)+p^2)}{x(s)(x(t)+p^2)}\\
&\qquad=\log\left\{e^{-2q(t-s)}\frac{(2q+p)^2e^{2qt}-p^2}{(2q+p)^2e^{2qs}-p^2}
\right\},
\end{split}
\end{equation*}
so that
\[
e^{-\int_{s}^{t}k(u,u)du}
=e^{-p(t-s)}e^{-2q(t-s)}\frac{(2q+p)^2e^{2qt}-p^2}{(2q+p)^2e^{2qs}-p^2}.
\]
Thus
\[
l(t,s)=(2q+p)pe^{-(p+q)(t-s)}
\frac{(2q+p)e^{2qs}-p}{(2q+p)^2e^{2qs}-p^2}
\]
or (\ref{eq:1.4}), as desired.
\end{proof}

\section{Filtering equations}

\label{sec:3}

\subsection{General result}

\label{sec:3.1}

Let $(X(t),U(t))_{0\le t\le T}$ be a two-dimensional centered continuous
Gaussian process. The process $X(\cdot)$ represents the state process, 
while 
$U(\cdot)$ is another process which is related to the dynamics of 
$X(\cdot)$.

Let $(B(t))_{0\leq t\leq T}$ be a one-dimensional Brownian motion that is
independent of $(X,U)$. We define the processes $V(\cdot )$ and $\alpha
(\cdot )$ by (\ref{eq:2.1}) and (\ref{eq:2.3}), respectively. In this
subsection, we assume that $l(t,s)$ in (\ref{eq:2.3}) is an arbitrary
Volterra-type bounded measurable function (i.e., $l(t,s)=0$ for $s>t$)
rather than the special form (\ref{eq:1.4}). Of course, the function 
$l(t,s)$
in (\ref{eq:1.4}) satisfies this assumption. We consider the observation
process $Y(\cdot )$ satisfying 
\begin{equation*}
Y(t)=\int_{0}^{t}\mu (s)X(s)ds+V(t)\qquad (0\leq t\leq T),
\end{equation*}
where $\mu (\cdot )$ is a bounded measurable deterministic function on 
$[0,T]
$ such that $\mu (t)\neq 0$ for $0\leq t\leq T$.

As in Section \ref{sec:2}, let $\mathcal{F}^Y$ be the augmented filtration
generated by $Y(\cdot)$. For $d$-dimensional column vector processes $%
(A(t))_{0\le t\le T}$ and $(C(t))_{0\le t\le T}$, we write 
\begin{gather*}
\widehat{A}(t)=\mathbb{E}[A(t)|\mathcal{F}^Y(t)] \qquad (0\le t\le T), \\
\Gamma_{AC}(t,s)=\mathbb{E}[A(t)C^{*}(s)] \qquad (0\le s\le t\le T),
\end{gather*}
where $*$ denotes the transposition. Notice that $\Gamma_{AC}(t,s)\in 
\mathbf{R}^{d\times d}$.

We put 
\begin{equation*}
Z(t)=(X(t),U(t),\alpha (t))^{\ast }\qquad (0\leq t\leq T),
\end{equation*}
and define the error matrix $P(t,s)\in \mathbf{R}^{3\times 3}$ by 
\begin{equation*}
P(t,s)=\mathbb{E}[Z(t)(Z(s)-\widehat{Z}(s))^{\ast }]\qquad (0\leq s\leq
t\leq T).
\end{equation*}

The next theorem gives an answer to the filtering problem for the partially
observable process $(X(t),Y(t))_{0\leq t\leq T}$. It turns out that this
will be useful in the filtering problem for (\ref{eq:1.5}) for example.

\begin{thm}
\label{thm:3.1} The filter $\widehat{Z}(\cdot)$ satisfies the stochastic
integral equation 
\begin{equation}
\widehat{Z}(t)=\int_0^t \{P(t,s)+D(t,s)\}a(s) \{dY(s)-a^{*}(s)\widehat{Z}%
(s)ds\},  \label{eq:3.1}
\end{equation}
and the error matrix $P(t,s)$ is the unique solution to the following 
matrix
Riccati-type integral equation such that $P(t,t)$ is symmetric and
nonnegative definite for $0\le t\le T$: 
\begin{equation}
\begin{split}
P(t,s)&=- \int_0^s \{P(t,r)+D(t,r)\}a(r)a^{*}(r)\{P(s,r)+D(s,r)\}^{*}dr \\
&\qquad\qquad\qquad\qquad\qquad\quad +\Gamma_{ZZ}(t,s)\qquad (0\le s\le 
t\le
T),
\end{split}
\label{eq:3.2}
\end{equation}
where 
\begin{equation*}
D(t,s)=\left(
\begin{matrix}
0 & 0 & 0\cr 0 & 0 & 0\cr l(t,s)/\mu(s) & 0 & 0
\end{matrix}
\right),\quad a(s)=\left(
\begin{matrix}
\mu(s)\cr 0\cr -1
\end{matrix}
\right).
\end{equation*}
\end{thm}

\begin{proof}
Since $(X,U)$ is independent of $B$, $(X,U,\alpha,Y)$ forms 
a Gaussian system. 
We have
\begin{equation*}
Y(t)=\int_{0}^{t}\{\mu(s)X(s)-\alpha(s)\}ds+B(t).
\end{equation*}
Thus we can define the innovation process $I(\cdot)$ by 
\begin{equation*}
 I(t)=Y(t)- \int_0^t(\mu(s)\widehat{X}(s)ds-\widehat{\alpha}(s)) 
 \qquad (0\le t\le T),
\end{equation*}
which is a Brownian motion satisfying $\mathcal{F}^Y=\mathcal{F}^I$ 
(cf.~\cite[Theorem 7.16]{LS}). 
Notice that $I(\cdot)$ can be written as 
\begin{equation}
 I(t)=\int_0^t (Z(s)-\widehat{Z}(s))^{*}a(s)ds + B(t).
\label{eq:3.3}
\end{equation}
By \cite[Corollary of Theorem 7.16]{LS}, there exists 
an $\mathbf{R}^3$-valued Volterra-type function 
$F(t,s)=(F_1(t,s),F_2(t,s),F_3(t,s))^{*}$ on $[0,T]^2$ such that 
\begin{align}
 &\int_0^t|F(t,s)|^2 ds < +\infty\quad (0\le t\le T), \nonumber \\
 &\widehat{Z}(t)=\int_0^t F(t,s)dI(s) \quad (0\le t\le T), 
\label{eq:3.4} 
\end{align}
where $|\cdot|$ denotes the Euclidean norm

Now let $g(t)=(g_1(t),g_2(t),g_3(t))$ be an arbitrary 
bounded measurable row-vector function on $[0,T]$. 
Then, for $t\in [0,T]$, 
\begin{equation*}
 \mathbb{E}\left[\int_0^t g(s)dI(s)\cdot(Z(t)-\widehat{Z}(t))\right]=0. 
\end{equation*}
From this, (\ref{eq:3.3}), (\ref{eq:3.4}) and the fact that $(X,U)$ and $B$ 
are independent, we have
\begin{equation*}
\begin{split}
 &\int_0^t g(s)F(t,s)ds
 = \mathbb{E}\left[\int_0^t g(s)dI(s)\cdot Z(t)\right] \\
 &= \mathbb{E}\left[\int_0^t g(s)\{(Z(s)-\widehat{Z}(s))^{*}a(s)ds 
 + dB(s)\}\cdot Z(t)\right]\\
 &= \int_0^t g(s)\mathbb{E}[Z(t)(Z(s)-\widehat{Z}(s))^{*}]a(s)ds 
   + \int_0^t g_3(s)l(t,s)ds \\
 &= \int_0^t g(s)P(t,s)a(s)ds 
   + \int_0^t g(s)D(t,s)a(s)ds.
\end{split}
\end{equation*}
Since $g(\cdot)$ is arbitrary, we deduce that $F(t,s)=(P(t,s)+D(t,s))a(s)$ 
or
\begin{equation}
\widehat{Z}(t)=\int_0^t (P(t,s)+D(t,s))a(s)dI(s) \quad (0\le t\le T).
\label{eq:3.5}
\end{equation}

The SDE (\ref{eq:3.1}) follows from (\ref{eq:3.5}) 
and the representation
\begin{equation*}
I(t)=Y(t)-\int_{0}^{t}a^{*}(s)\widehat{Z}(s)ds\quad (0\le t\le T).
\end{equation*}
The equation (\ref{eq:3.2}) follows from (\ref{eq:3.5}) and the equality
\begin{equation*}
P(t,s)=\mathbb{E}[Z(t)Z^{*}(s)]-
   \mathbb{E}[\widehat{Z}(t)\widehat{Z}^{*}(s)].
\end{equation*}
The matrix $P(t,t)$ is clearly symmetric and nonnegative definite. 
Finally, the uniqueness of the solution to (\ref{eq:3.2}) follows 
from Proposition \ref{prop:3.2} below.  
\end{proof}

\begin{prop}
\label{prop:3.2} The solution $P(t,s)$ to the matrix integral equation 
$(\ref
{eq:3.2})$ such that $P(t,t)$ is symmetric and nonnegative definite for $%
0\le t\le T$ is unique.
\end{prop}

\begin{proof}
By coninuity, there exists a positive 
constant $C(T)$ such that 
\begin{equation*}
\|\Gamma_{ZZ}(t,s)\|\le C(T) \qquad (0\le s\le t\le T),  
\end{equation*}
where $\|A\|:=\left(\mbox{trace}(A^{*}A)\right)^{1/2}$ 
for $A\in\mathbf{R}^{3\times 3}$. 
Let $P$ be a solution to (\ref{eq:3.2}) such that 
$P(t,t)$ is symmetric and nonnegative definite for $t\in [0,T]$. 
We put $Q(t,s)=P(t,s)+D(t,s)$. 
Then (\ref{eq:3.2}) with $s=t$ is
\begin{equation*}
 \Gamma_{ZZ}(t,t)- P(t,t)= \int_0^t Q(t,r)a(r)a^{*}(r)Q^{*}(t,r)dr. 
\end{equation*}
From this, we have
\begin{align*}
 &\int_0^t |Q(t,r)a(r)|^2 dr 
 = \int_0^t \mbox{trace}
  \left\{Q(t,r)a(r)a^{*}(r)Q^{*}(t,r)\right\}dr \\ 
 &\qquad=\mbox{trace}(\Gamma_{ZZ}(t,t)-P(t,t))
 \le \mbox{trace}(\Gamma_{ZZ}(t,t))\\
 &\qquad\le \sqrt{3}\|\Gamma_{ZZ}(t,t)\|
 \le \sqrt{3}C(T). 
\end{align*}
Therefore, $\|P(t,s)\|$ is at most 
\begin{align*}
 &\|\Gamma_{ZZ}(t,s)\| + 
  \int_0^s \|Q(t,r)a(r)a^{*}(r)Q^{*}(s,r)\|dr \\
 &\le \|\Gamma_{ZZ}(t,s)\| + 
  \left(\int_0^t |Q(t,r)a(r)|^2dr\right)^{1/2}\cdot 
  \left(\int_0^s |Q(s,r)a(r)|^2dr\right)^{1/2} \\
 &\le (1+\sqrt{3})C(T).
\end{align*}

Let $P_1$ and $P_2$ be two solutions of (\ref{eq:3.2}). We define 
$Q_i(t,s)=P_i(t,s)+D(t,s)$ for $i=1,2$.  
We put $P_i(t,s)=0$ for $s>t$ and $i=1,2$. 
Since $\mu$ and $l$ are bounded, it follows from the above estimate that
there exists a positive constant $K(T)$ satisfying
\begin{equation*}
 \|a(s)a^{*}(s)Q_i(t,s)\|\le K(T) \qquad (0\le s\le t\le T,\ i=1,2). 
\end{equation*}
It follows that
\begin{equation*}
\begin{split}
&\|Q_1(t,r)a(r)a^{*}(r)Q_1(s,r)^{*}-Q_2(t,r)a(r)a^{*}(r)Q_2(s,r)^{*}\|\\
&\qquad \le \|Q_1(t,r)a(r)a^{*}(r)(Q_1(s,r)^{*}-Q_2(s,r)^{*})\|\\
&\qquad\qquad\qquad+\|(Q_1(t,r)-Q_2(t,r))a(r)a^{*}(r)Q_2(s,r)^{*}\|\\
&\qquad\le 2K(T)\|Q_1(s,r)-Q_2(t,r)\|
=2K(T)\|P_1(s,r)-P_2(t,r)\|.
\end{split}
\end{equation*}
From this and (\ref{eq:3.2}), we obtain
\begin{equation*}
 \sup_{0\le t\le T}\|P_1(t,s)-P_2(t,s)\|\le 
 2K(T)\int_0^s\sup_{0\le t\le T}\|P_1(t,r)-P_2(t,r)\|dr. 
\end{equation*}
Therefore, Gronwall's lemma gives 
\begin{equation*}
 \sup_{0\le t\le T}\|P_1(t,s)-P_2(t,s)\|=0 \qquad (0\le s\le T). 
\end{equation*}
Thus the uniqueness follows.
\end{proof}

\begin{rem}
We consider the case in which $\alpha =0$ and the state process $X(\cdot )$
is the Ornstein-Uhlenbeck process satisfying 
\begin{equation*}
dX(t)=\theta X(t)dt+\sigma dW(t),\quad X(0)=0,
\end{equation*}
where $\theta ,\sigma \neq 0$ and $W(\cdot )$ is a one-dimensional Brownian
motion that is independent of $B(\cdot )$. We also assume that $\mu (\cdot
)=\mu ,$ a constant. Then $X(t)=\sigma \int_{0}^{t}e^{\theta (t-u)}dW(u)$
and 
\begin{equation*}
\Gamma _{XX}(t,s)=\frac{\sigma ^{2}}{2\theta }(e^{\theta (t+s)}-e^{\theta
(t-s)})\quad (0\leq s\leq t\leq T).
\end{equation*}
By Theorem \ref{thm:3.1}, we have 
\begin{gather}
\widehat{X}(t)=\int_{0}^{t}\mu P_{XX}(t,s)(dY(s)-\mu \widehat{X}(s)ds),
\label{eq:3.6} \\
P_{XX}(t,s)=\Gamma _{XX}(t,s)-\int_{0}^{s}\mu ^{2}P_{XX}(t,r)P_{XX}(s,r)dr,
\label{eq:3.7}
\end{gather}
where $P_{XX}(t,s)=\mathbb{E}[X(t)(X(s)-\widehat{X}(s))]$ for $0\leq s\leq
t\leq T$. Let $\mathcal{F}(t)$, $0\leq t\leq T$, be the $\mathbb{P}$%
-augmentation of the filtration generated by $(W(t),B(t))_{0\leq t\leq T}$.
Then $P_{XX}(t,s)$ is 
\begin{equation*}
\begin{split}
\mathbb{E}[X(t)(X(s)-\widehat{X}(s))]& 
=\mathbb{E}[\mathbb{E}[X(t)|\mathcal{F%
}(s)](X(s)-\widehat{X}(s))] \\
& =\mathbb{E}[e^{\theta (t-s)}X(s)(X(s)-\widehat{X}(s))]=e^{\theta
(t-s)}\gamma (s)
\end{split}
\end{equation*}
with $\gamma (s)=\mathbb{E}[X(s)(X(s)-\widehat{X}(s))]$. Thus 
(\ref{eq:3.6})
is reduced to 
\begin{equation}
d\widehat{X}(t)=(\theta -\mu ^{2}\gamma (t))\widehat{X}(t)dt+\mu \gamma
(t)dY(t).  \label{eq:3.8}
\end{equation}
Differentiating (\ref{eq:3.7}) in $s=t$, we get 
\begin{equation}
\frac{d\gamma (t)}{dt}=\sigma ^{2}+2\theta \gamma (t)-\mu ^{2}\gamma 
(t)^{2}.
\label{eq:3.9}
\end{equation}
The equations (\ref{eq:3.8}) and (\ref{eq:3.9}) are the well-known
Kalman-Bucy filtering equations (see \cite{KB}, \cite{BJ}, \cite{D}, 
\cite{J}
and \cite{LS}).
\end{rem}

\subsection{Linear systems with memory}

\label{sec:3.2}

Let $(X(t),Y(t))_{0\le t\le T}$ be the two-dimensional partially observable
process satisfying (\ref{eq:1.5}). For $j=1,2$, the noise process 
$V_j(\cdot)
$ there is described by (\ref{eq:1.1}) with $W(\cdot)=W_j(\cdot)$ and $%
(p,q)=(p_j,q_j)$ satisfying (\ref{eq:1.2}). The Brownian motions 
$W_1(\cdot)$
and $W_2(\cdot)$, whence $V_1(\cdot)$ and $V_2(\cdot)$, are independent. In 
(%
\ref{eq:1.5}), the coefficients $\theta, \sigma, \mu\in\mathbf{R}$ with $%
\mu\neq 0$ are known constants and the initial value $X_0$ is a Gaussian
random variable that is independent of $(V_1,V_2)$. The processes 
$X(\cdot)$
and $Y(\cdot)$ represent the state and the observation, respectively.

By Theorem \ref{thm:2.1}, we have $V_j(t)= B_j(t)- \int_0^t \alpha_j(s)ds$
with $\alpha_j(t)=\int_{0}^{t}l_j(t,s)dB_j(s)$ for $j=1,2$ and $0\le t\le 
T$%
, where $B_1(\cdot)$ and $B_2(\cdot)$ are two independent Brownian motions,
and, for $0\le s\le t\le T$, 
\begin{equation*}
l_j(t,s)= p_je^{-(p_j+q_j)(t-s)}\left\{1-\frac{2p_jq_j} {%
(2q_j+p_j)^2e^{2q_js}-p_j^2}\right\}.
\end{equation*}
We put $l_j(t)=l_j(t,t)$ for $j=1,2$, that is, 
\begin{equation}  \label{eq:3.10}
l_j(t)= p_j\left\{1-\frac{2p_jq_j}{(2q_j+p_j)^2e^{2q_js}-p_j^2}\right\}
\quad (j=1,2,\ 0\le t\le T).
\end{equation}
It holds that $l_j(t,s)=e^{-r_j(t-s)}l_j(s)$, where $r_j=p_j+q_j$.

We denote by $\mathcal{F}(t)$, $0\le t\le T$, the $\mathbb{P}$-augmentation
of the filtration $\sigma(X_0, (V_1(s),V_2(s))_{0\le s\le t})$, $0\le t\le 
T$%
.

\begin{lem}
\label{lem:3.3} For $0\le s\le t\le T$, we have 
\begin{equation*}
\mathbb{E}[X(t)|\mathcal{F}(s)] =e^{\theta(t-s)}X(s)-\sigma
b(t-s)\alpha_1(s),
\end{equation*}
where 
\begin{equation*}
b(t)= 
\begin{cases}
(e^{\theta t}- e^{-r_1 t})/(\theta + r_1) & (\theta + r_1 \neq 0), \\ 
te^{\theta t} & (\theta + r_1 = 0).
\end{cases}
\end{equation*}
\end{lem}

\begin{proof}
 For $t\in [0,T]$, $\int_0^t e^{\theta (t-s)}\alpha_1(s)ds$ is 
\begin{equation*}
 \int_0^tl_1(u)\left\{\int_u^t e^{\theta(t-s)}e^{-r_1(s-u)}ds
 \right\}dB_1(u)= \int_0^t b(t-u)l_1(u)dB_1(u). 
\end{equation*}
Since $X(t)=e^{\theta t}X_0 + \sigma\int_0^t e^{\theta(t-u)}dV_1(u)$ 
or
\begin{equation*}
 X(t)=e^{\theta t}X_0 + \sigma\int_0^t e^{\theta(t-u)}dB_1(u)
      -\sigma\int_0^t e^{\theta(t-u)}\alpha_1(u)du, 
\end{equation*}
$\mathbb{E}[X(t)|\mathcal{F}(s)]$ with $s\le t$ is equal to
\begin{align*}
 &e^{\theta t}X_0+ 
   \sigma e^{\theta(t-s)}\int_0^s e^{\theta(s-u)}dB_1(u)
   -\sigma \int_0^s b(t-u)l_1(u)dB_1(u) \\
 &= e^{\theta(t-s)}X(s) -\sigma\int_0^s 
     \left\{b(t-u)-e^{\theta(t-s)}b(s-u)\right\}l_1(u)dB_1(u). 
\end{align*}
However, by elementary calculation, we have
\begin{equation*}
 b(t-u)-e^{\theta(t-s)}b(s-u)= b(t-s)e^{-r_1(s-u)} 
 \qquad (0\le u\le s\le t).   
\end{equation*}
Thus the lemma follows.
\end{proof}

We put, for $0\le t\le T$, 
\begin{gather*}
F=\left(
\begin{matrix}
-\theta & \sigma & 0\cr 0 & r_1 & 0\cr 0 & 0 & r_2
\end{matrix}
\right), \quad D(t)=\left(
\begin{matrix}
0 & 0 & 0\cr 0 & 0 & 0\cr \mu^{-1}l_2(t) & 0 & 0
\end{matrix}
\right), \quad a=\left(
\begin{matrix}
\mu\cr 0\cr -1
\end{matrix}
\right), \\
G(t)=\left(
\begin{matrix}
\sigma^2 & \sigma l_1(t) & 0\cr \sigma l_1(t) & l_1^2(t) & 0\cr 0 & 0 & 0
\end{matrix}
\right), \quad H(t)=\left(
\begin{matrix}
-\theta & \sigma & 0\cr 0 & r_1 & 0\cr \mu l_2(t) & 0 & r_2-l_2(t)
\end{matrix}
\right).
\end{gather*}
We also put 
\begin{equation*}
Z(t)=(X(t),\alpha_1(t),\alpha_2(t))^{*}\qquad (0\le t\le T).
\end{equation*}
Recall that $\Gamma_{ZZ}(0)=\mathbb{E}[Z(0)Z^{*}(0)]$ and that 
\begin{equation*}
\widehat{Z}(t)=\mathbb{E}[Z(t)|\mathcal{F}^Y(t)] \qquad (0\le t\le T).
\end{equation*}
We define the error matrix $P(t)\in\mathbf{R}^{3\times 3}$ by 
\begin{equation*}
P(t)=\mathbb{E}[Z(t)(Z(t)-\widehat{Z}(t))^{*}]\qquad (0\le t\le T).
\end{equation*}

Here is the solution to the optimal filtering problem for (\ref{eq:1.5}).

\begin{thm}
\label{thm:3.4} The filter $\widehat{Z}(\cdot )$ satisfies the stochastic
differential equation 
\begin{equation}
\begin{split}
d\widehat{Z}(t)& =-\{F+(P(t)+D(t))aa^{\ast }\}\widehat{Z}(t)dt \\
& \qquad \qquad \qquad \qquad +(P(t)+D(t))adY(t)\quad (0\leq t\leq T),
\end{split}
\label{eq:3.11}
\end{equation}
with $\widehat{Z}(0)=(\mathbb{E}[X_{0}],0,0)^{\ast }$, and $P(\cdot )$
follows the matrix Riccati equation 
\begin{equation}
\begin{split}
\frac{dP(t)}{dt}& =G(t)-H(t)P(t)-P(t)H(t)^{\ast }-P(t)aa^{\ast }P(t) \\
& \qquad \qquad \qquad \qquad \qquad \qquad \qquad \qquad (0\leq t\leq T)
\end{split}
\label{eq:3.12}
\end{equation}
with $P_{ij}(0)=\delta _{i1}\delta _{j1}\mathbb{E}[(X_{0})^{2}]$ for $%
i,j=1,2,3$.
\end{thm}

\begin{proof}
For $0\le s\le t\le T$, we put 
\begin{equation*}
 P(t,s)=\mathbb{E}[Z(t)(Z(s)-\widehat{Z}(s))^{*}].
 \end{equation*}
Then we have $P(t)=P(t,t)$. We also put, for $0\le s\le t\le T$,
\begin{align*}
&D(t,s)=e^{-r_2(t-s)}D(s),\\
&Q(t,s)=P(t,s)+D(t,s),\quad Q(s)=Q(s,s)=P(s)+D(s).
\end{align*}

By
\begin{equation*}
P(t,s)=\mathbb{E}[\mathbb{E}[Z(t)|\mathcal{F}(s)](Z(s)-\widehat{Z}(s))^{*}]
\end{equation*}
and Lemma \ref{lem:3.3}, we have $P(t,s)=M(t-s)P(s)$, 
where
\begin{equation*}
 M(t)=\left(\begin{matrix}
	     e^{\theta t} & -\sigma b(t) & 0\cr
	     0 & e^{-r_1 t} & 0\cr
	     0 & 0 & e^{-r_2 t}
	    \end{matrix}\right)
\end{equation*}
with $b(\cdot)$ as in Lemma \ref{lem:3.3}. 
We also see that $D(t,s)=M(t-s)D(s)$. Combining, 
$Q(t,s)= M(t-s)Q(s)$. However, $M(t)=e^{-tF}$ since 
$dM(t)/dt=-FM(t)$ and $M(0)$ is the unit matrix. Thus we obtain
\begin{equation}
Q(t,s)=e^{-(t-s)F}Q(s).
\label{eq:3.13}
\end{equation}

From (\ref{eq:3.13}) and Theorem \ref{thm:3.1} with 
$U=\alpha_1$ and $\alpha=\alpha_2$, it follows that
\begin{gather}
\widehat{Z}(t)= \int_0^t e^{-(t-s)F}Q(s)a\{dY(s)-a^{*}\widehat{Z}(s)ds\}, 
\label{eq:3.14}\\
P(t)=\Gamma_{ZZ}(t)-\int_0^t e^{-(t-u)F}Q(u)aa^{*}
   Q^{*}(u)e^{-(t-u)F^{*}}du.
\label{eq:3.15}
\end{gather}
The SDE (\ref{eq:3.11}) follows easily from (\ref{eq:3.14}).

Differentiating both sides of (\ref{eq:3.15}) 
yields
\begin{equation}
\begin{split}
 \dot{P}&=\dot{\Gamma}_{ZZ}+F(\Gamma_{ZZ}-P)+ (\Gamma_{ZZ}-P)F^{*}
  -Qaa^{*}Q^{*} \\
 &=\dot{\Gamma}_{ZZ}+F\Gamma_{ZZ}+\Gamma_{ZZ}F^{*}-Daa^{*}D
   -HP-PH^{*}- Paa^{*}P. 
\end{split}
\label{eq:3.16}
\end{equation}
Since
\begin{equation*}
 dZ(t)=-FZ(t)dt+ dR(t)
\end{equation*}
with
\begin{equation*}
 R(t)=\left(\sigma B_1(t), \int_0^t l_1(s)dB_1(s), 
       \int_0^t l_2(s)dB_2(s)\right)^{*},
\end{equation*}
we see by integration by parts that $Z(t)Z(t)^{*}-Z(0)Z(0)^{*}$ is equal to
\begin{equation*}
\int_{0}^{t}Z(s)dZ(s)^{*}+\int_{0}^{t}dZ(s)Z(s)^{*}+\mathbb{E}[R(t)R^{*}(t)].
\end{equation*}
It follows that $\Gamma_{ZZ}(t)-\Gamma_{ZZ}(0)$ is equal to
\begin{equation*}
\begin{split}
 &\mathbb{E}\left[\int_{0}^{t}Z(s)dZ(s)^{*}\right]
+\mathbb{E}\left[\int_{0}^{t}dZ(s)Z(s)^{*}\right]+\mathbb{E}[R(t)R^{*}(t)]\\
 &\qquad=-\int_0^t\Gamma_{ZZ}(s)F^{*}ds-\int_0^t F\Gamma_{ZZ}(s)ds
 +\mathbb{E}[R(t)R^{*}(t)].
\end{split}
\end{equation*}
Thus
\begin{equation*}
\begin{split}
 &\quad\Gamma_{ZZ}(t)+\int_0^t F\Gamma_{ZZ}(s)ds 
  +\int_0^t\Gamma_{ZZ}(s)F^{*}ds - \Gamma_{ZZ}(0)\\ 
 &\qquad\quad=\mathbb{E}[R(t)R^{*}(t)]
  =\int_0^t\left(G(s)+D(s)aa^{*}D^{*}(s)\right)ds. 
\end{split}
\end{equation*}
This and (\ref{eq:3.16}) yield (\ref{eq:3.12}).
\end{proof}

\begin{rem}
We can easily extend Theorem \ref{thm:3.4} to a more general setting where
the functions $l_{j}(t,s)$ take the form $l_{j}(t,s)=e^{c(t-s)}g(s)$.
\end{rem}

\begin{cor}
\label{cor:3.5} We assume that $p_2=0$, i.e., $V_2(\cdot)=W_2(\cdot)$. Let 
$%
Z(t)=(X(t),\alpha_1(t))^{*}$ and $P(t)=\mathbb{E}[Z(t)(Z(t)-\widehat{Z}%
(t))^{*}]\in\mathbf{R}^{2\times 2}$. Then the filter $\widehat{Z}(\cdot)$
and the error matrix $P(\cdot)$ satisfy, respectively, 
\begin{gather*}
d\widehat{Z}(t)=-\{F+ P(t)aa^{*}\}\widehat{Z}(t)dt + P(t)adY(t),\quad 
\widehat{Z}(0)=(\mathbb{E}[X_0], 0)^{*}, \\
\frac{dP(t)}{dt}=G(t)-FP(t)-P(t)F^{*}-P(t)aa^{*}P(t), \quad
P_{ij}(0)=\delta_{i1}\delta_{j1}\mathbb{E}[(X_0)^2],
\end{gather*}
where 
\begin{equation*}
F=\left(
\begin{matrix}
-\theta & \sigma\cr 0 & r_1
\end{matrix}
\right), \quad a=\left(
\begin{matrix}
\mu\cr 0
\end{matrix}
\right), \quad G(t)=\left(
\begin{matrix}
\sigma^2 & \sigma l_1(t)\cr \sigma l_1(t) & l_1(t)^2
\end{matrix}
\right).
\end{equation*}
\end{cor}

\begin{cor}
\label{cor:3.6} We assume that $p_1=0$, i.e., $V_1(\cdot)=W_1(\cdot)$. Let 
$%
Z(t)=(X(t),\alpha_2(t))^{*}$ and $P(t)=\mathbb{E}[Z(t)(Z(t)-\widehat{Z}%
(t))^{*}] \in\mathbf{R}^{2\times 2}$. Then the filter $\widehat{Z}(\cdot)$
and the error matrix $P(\cdot)$ satisfy, respectively, 
\begin{gather*}
d\widehat{Z}(t)=-\{F+(P(t)+D(t))aa^{*}\}\widehat{Z}(t)dt + 
(P(t)+D(t))adY(t),
\\
\frac{dP(t)}{dt}=G-H(t)P(t)-P(t)H(t)^{*}-P(t)aa^{*}P(t), \quad ,
\end{gather*}
with $\widehat{Z}(0)=(\mathbb{E}[X_0], 0)^{*}$ and $P_{ij}(0)=\delta_{i1}%
\delta_{j1}\mathbb{E}[(X_0)^2]$, where 
\begin{gather*}
F=\left(
\begin{matrix}
-\theta & 0\cr 0 & r_2
\end{matrix}
\right), \quad D(t)=\left(
\begin{matrix}
0 & 0\cr \mu^{-1}l_2(t) & 0
\end{matrix}
\right), \quad a=\left(
\begin{matrix}
\mu\cr -1
\end{matrix}
\right), \\
G=\left(
\begin{matrix}
\sigma^2 & 0\cr 0 & 0
\end{matrix}
\right), \quad H(t)=\left(
\begin{matrix}
-\theta & 0\cr \mu l_2(t) & r_2-l_2(t)
\end{matrix}
\right).
\end{gather*}
\end{cor}

\begin{ex}
We consider the estimation problem of the value of a signal $\rho $ from 
the
observation process $(Y(t))_{0\leq t\leq T}$ given by 
\begin{equation*}
dY(t)=\rho dt+dV(t),\quad Y(0)=0,
\end{equation*}
where $V(\cdot )$ and $\alpha (\cdot )$ are as in Section \ref{sec:2}. We
assume that $\rho $ is a random variable having the normal distribution 
with
mean zero and variance $v^{2}$. This is the special case $\theta =\sigma 
=0$
of the setting of Corollary \ref{cor:3.6}. Let $r=p+q$ and $l(\cdot )$ be 
as
above. Let $H(t)$ and $a$ be as in Corollary \ref{cor:3.6} with $\mu =1$ 
and 
$\theta =0$. We define $P(t)=(P_{ij}(t))_{1\leq i,j\leq 2}$ by 
$P(t)=\mathbb{%
E}[Z^{\ast }(t)(Z(t)-\widehat{Z}(t))]$ with $Z(t)=(\rho ,\alpha (t))^{\ast 
}$%
. Then, by Corollary \ref{cor:3.6}, the filter $(\widehat{\rho 
}(t),\widehat{%
\alpha }(t))$ satisfies 
\begin{align*}
d\widehat{\rho }(t)& =-\{P_{11}(t)-P_{12}(t)\}\{\widehat{\rho 
}(t)-\widehat{%
\alpha }(t)\}dt+\{P_{11}(t)-P_{12}(t)\}dY(t), \\
d\widehat{\alpha }(t)& =[\{P_{21}(t)-P_{22}(t)+l(t)-r\}\widehat{\alpha }%
(t)-\{P_{21}(t)-P_{22}(t)+l(t)\}\widehat{\rho }(t)]dt \\
& \qquad \qquad \qquad \qquad \qquad +\{P_{21}(t)-P_{22}(t)+l(t)\}dY(t)
\end{align*}
with $(\widehat{\rho }(0),\widehat{\alpha }(0))=(\mathbb{E}[X_{0}],0)$, and
the error matrix $P(\cdot )$ follows 
\begin{equation*}
\frac{dP(t)}{dt}=-H(t)P(t)-P(t)H(t)^{\ast }-P(t)aa^{\ast }P(t),\quad
P_{ij}(0)=\delta _{i1}\delta _{j1}\mathbb{E}[\rho ^{2}].
\end{equation*}
By the linearization method as described in \cite[Chapter 5]{BJ}, we can
solve the equation for $P(\cdot )$ to obtain 
\begin{equation*}
P(t)=\frac{v^{2}}{1+v^{2}\eta (t)+v^{2}\xi (t)\phi (t)}\left( 
\begin{matrix}
1 & -\phi (t)/\psi (t) \\ 
-\phi (t)/\psi (t) & v^{2}\phi (t)^{2}/\psi (t)^{2}
\end{matrix}
\right) ,
\end{equation*}
where 
\begin{gather*}
\psi (t)=\exp \left\{ \int_{0}^{t}(r-l(s))ds\right\} ,\quad \phi
(t)=\int_{0}^{t}l(s)\psi (s)ds, \\
\xi (t)=\int_{0}^{t}\frac{\psi (s)+\phi (s)}{\psi (s)^{2}}ds, \\
\eta (t)=\int_{0}^{t}\left\{ 1-\psi (s)\xi (s)+\frac{\phi (s)}{\psi (s)}%
\right\} ds.
\end{gather*}
The analytical forms of $\psi $, $\phi $, $\xi $ and $\eta $ can be 
derived.
We omit the details.
\end{ex}

\section{Application to finance}

\label{sec:4}

In this section, we apply the results in the previous section to the 
problem
of expected utility maximization for an investor with partial observations.

Let $(V_j(t))_{0\le t\le T}$, $(\alpha_j(t))_{0\le t\le T}$, $j=1,2$, be as
in Section \ref{sec:3}. In particular, $V_1(\cdot)$ and $V_2(\cdot)$ are
independent. We consider the financial market consisting of a share of the
money market with price $S_0(t)$ at time $t\in [0,T]$ and a stock with 
price 
$S(t)$ at time $t\in [0,T]$. The stock price $S(\cdot)$ is governed by the
stochastic differential equation 
\begin{equation}  \label{eq:4.1}
dS(t)=S(t)\{U(t)dt + \eta dV_2(t)\}, \quad S(0)=s_0,
\end{equation}
where $s_0$ and $\eta$ are positive constants and $U(\cdot)$ is a Gaussian
process following 
\begin{equation}  \label{eq:4.2}
dU(t)= \{\delta +\theta U(t)\}dt + \sigma dV_1(t), \quad U(0)=\rho.
\end{equation}
The parameters $\theta$, $\delta$ and $\sigma$ are constants, and $\rho$ is
a Gaussian random variable that is independent of $(V_1,V_2)$. For
simplicity, we assume that 
\begin{equation*}
S_0(\cdot)=1,\quad \eta=1,\quad \delta=0.
\end{equation*}

Let $\mathcal{F}(t)$, $0\le t\le T$, be the $\mathbb{P}$-augmentation of 
the
filtration generated by $(V_1(s), V_2(s))_{0\le s\le t}$ and $\rho$. Then 
$%
U(\cdot)$ is $\mathcal{F}$-adapted but not $\mathcal{F}^S$-adapted; recall
from Section \ref{sec:2} that $\mathcal{F}^S$ is the augmented filtration
generated by $S(\cdot)$. Suppose that we can observe neither the 
disturbance
process $V_2(\cdot)$ nor the drift process $U(\cdot)$ but only the price
process $S(\cdot)$. Thus only $\mathcal{F}^S$-adapted processes are
observable.

In many references such as \cite{DE}, \cite{DF}, and \cite{G}, the 
partially
observable model described by (\ref{eq:4.1}) and (\ref{eq:4.2}) with 
$V_j$'s
replaced by Brownian motions, i.e., $V_j=B_j$, is studied.

We consider the following expected logarithmic utility maximization from
terminal wealth: for given initial capital $x\in (0,\infty)$, 
\begin{equation}  \label{eq:4.3}
\mbox{maximize }\mathbb{E}[\log(X^{x,\pi}(T))] \quad \mbox{over all 
}\pi\in%
\mathcal{A}(x),
\end{equation}
where 
\begin{equation*}
\mathcal{A}(x)=\left\{ (\pi(t))_{0\le t\le T}: 
\begin{matrix}
\mbox{$\pi(\cdot)$ is $\mathbf{R}$-valued, 
   $\mathcal{F}^S$-progressively measurable,} \\ 
\int_0^T \pi^2(t)dt <\infty, \; X^{x,\pi}(t)\ge 0\;\; (0\le t\le T)\mbox{
a.s.}
\end{matrix}
\right\},
\end{equation*}
and 
\begin{equation}  \label{eq:4.4}
X^{x,\pi}(t)= x+ \int_0^t\frac{\pi(u)}{S(u)}dS(u).
\end{equation}
The value $\pi(t)$ is the dollar amount invested in the stock at time $t$,
whence $\pi(t)/S(t)$ is the number of units of stock held at time $t$. The
process $X^{x,\pi}(\cdot)$ is the wealth process associated with the
self-financing portfolio determined uniquely by $\pi(\cdot)$.

An analogue of the problem (\ref{eq:4.3}) for full observations is solved in 
\cite{AIK}. For related work, see \cite{KZ}, \cite{L1}, \cite{L2} and the
references therein. We solve the problem (\ref{eq:4.3}) by combining the
results above on filtering and the martingale method as described in \cite
{KS}.

Consider the process $(Y(t))_{0\le t\le T}$ defined by 
\begin{equation*}
Y(t)= \int_0^t U(s)ds + V_2(t) = B_2(t)+ \int_0^t (U(s)-\alpha_2(s))ds,
\end{equation*}
which we regard as the observation process. Let $(I(t))_{0\le t\le T}$ be
the innovation process associated with $Y(\cdot)$ that is given by 
\begin{equation*}
I(t)= Y(t)- \int_0^t (\widehat{U}(s)-\widehat{\alpha}_2(s))ds,
\end{equation*}
where $\widehat{A}(t)=\mathbb{E}[A(t)|\mathcal{F}^Y(t)]$ as in the previous
sections. The innovation process $I(\cdot)$ is a $\mathcal{F}^Y$-Brownian
motion satisfying $\mathcal{F}^S=\mathcal{F}^Y=\mathcal{F}^I$.

Let $(L(t))_{0\le t\le T}$ be the exponential $\mathcal{F}$-martingale 
given
by 
\begin{align*}
L(t)&=\exp\left\{-\int_0^t(U(s)-\alpha_2(s))dB_2(s)-\frac{1}{2}
\int_0^t(U(s)-\alpha_2(s))^2ds\right\} \\
&=\exp\left\{-\int_0^t(U(s)-\alpha_2(s))dY(s)+\frac{1}{2} \int_0^t(U(s)-%
\alpha_2(s))^2ds\right\}.
\end{align*}
By \cite[Chapter 7]{LS}, we find that, for $t\in [0,T]$, 
\begin{align*}
\widehat{L}(t)&=\exp\left\{-\int_0^t(\widehat{U}(s)- \widehat{\alpha}%
_2(s))dY(s)+\frac{1}{2} \int_0^t(\widehat{U}(s)-\widehat{\alpha}%
_2(s))^2ds\right\} \\
&=\exp\left\{-\int_0^t(\widehat{U}(s)-\widehat{\alpha}_2(s))dI(s) 
-\frac{1}{2%
} \int_0^t(\widehat{U}(s)-\widehat{\alpha}_2(s))^2ds\right\}.
\end{align*}
The process $(\widehat{L}(t))_{0\le t\le T}$ is a 
$\mathcal{F}^Y$-martingale.

For $x\in (0,\infty)$ and $\pi\in\mathcal{A}(x)$, we see from It\^{o}
formula that the process $(\widehat{L}(t)X^{x,\pi}(t))$ is a local 
$\mathcal{%
F}^Y$-martingale, whence a $\mathcal{F}^Y$-supermartingale since $%
X^{x,\pi}(\cdot)$ is nonnegative. It follows that, for $x\in (0,\infty)$, 
$%
\pi\in\mathcal{A}(x)$, and $y\in (0,\infty)$, 
\begin{align}  \label{eq:4.5}
\mathbb{E}[\log(X^{x,\pi}(T))]&\le 
\mathbb{E}[\log(X^{x,\pi}(T))-y\widehat{L}%
(T)X^{x,\pi}(T)]+ yx \\
&\le \mathbb{E}[\log\{1/(y\widehat{L}(T))\}-1]+ yx,  \notag
\end{align}
where, in the second inequality, we have used 
\begin{equation*}
\log(z)-yz\le \log(1/y)-1\quad (y,z\in (0,\infty)).
\end{equation*}
The equalities in (\ref{eq:4.5}) hold if and only if 
\begin{equation}  \label{eq:4.6}
X^{x,\pi}(T)=x/\widehat{L}(T)\quad \mbox{a.s.}
\end{equation}
Thus the portfolio process $\pi(\cdot)$ satisfying (\ref{eq:4.6}) is 
optimal.

Put 
\begin{equation*}
\pi_0(t)=x(\widehat{U}(t)-\widehat{\alpha}_2(t))/ \widehat{L}(t) \qquad
(0\le t\le T).
\end{equation*}
Since $(x/\widehat{L})(0)=x$ and 
\begin{equation*}
d(x/\widehat{L})(t) 
=\frac{x(\widehat{U}(t)-\widehat{\alpha}_2(t))}{\widehat{%
L}(t)}dY(t) =\frac{\pi_0(t)}{S(t)}dS(t),
\end{equation*}
we see from (\ref{eq:4.4}) that the process $\pi_0(\cdot)$ satisfies (\ref
{eq:4.6}), whence it is the desired optimal optimal portfolio process. It
also satisfies 
\begin{equation*}
\frac{\pi_0(t)}{X^{x,\pi_0}(t)}= \widehat{U}(t)-\widehat{\alpha}_2(t) 
\qquad
(0\le t\le T).
\end{equation*}

We put 
\begin{equation*}
a=(1,0,-1)^{*}, \quad Z(t)=(U(t),\alpha_1(t),\alpha_2(t))^{*}\quad (t\in
[0,T]).
\end{equation*}
We define the error matrix $P(t)\in\mathbf{R}^{3\times 3}$ by $\mathbb{E}%
[Z(t)(Z(t)-\widehat{Z}(t))^{*}]$. Combining the results above with Theorem 
\ref{thm:3.4} which describes the dynamics of $\widehat{U}$ and $\widehat{%
\alpha}_2$, we obtain the next theorem.

\begin{thm}
\label{thm:4.1} The optimal portfolio process $\pi_0$ for the problem 
$(\ref
{eq:4.3})$ is given by 
\begin{equation*}
\pi_0(t)=xa^{*}\widehat{Z}(t)/\widehat{L}(t) \qquad (0\le t\le T)
\end{equation*}
and satisfies 
\begin{equation*}
X^{x,\pi_0}(T)= x/\widehat{L}(T), \quad \frac{\pi_0(t)}{X^{x,\pi_0}(t)}=
a^{*}\widehat{Z}(t) \qquad (0\le t\le T).
\end{equation*}
The filter $\widehat{Z}(\cdot)$ follows 
\begin{equation*}
d\widehat{Z}(t)=-\{F+(P(t)+D(t))aa^{*}\}\widehat{Z}(t)dt + 
(P(t)+D(t))adY(t)
\end{equation*}
with $\widehat{Z}(0)=(\mathbb{E}[\rho], 0, 0)^{*}$, and the error matrix $%
P(\cdot)$ satisfies the matrix Riccati equation 
\begin{equation*}
\frac{dP(t)}{dt}=G(t)-H(t)P(t)-P(t)H(t)^{*}-P(t)aa^{*}P(t)
\end{equation*}
with $P_{ij}(0)=\delta_{i1}\delta_{j1}\mathbb{E}[\rho^2]\ (i,j=1,2,3)$,
where $F$, $D(t)$, $G(t)$ and $H(t)$ are as in Theorem\/ $\ref{thm:3.4}$
with $\mu=1$.
\end{thm}

\section{Parameter estimation}

\label{sec:5}

Let $V(\cdot )$ be the process given by (\ref{eq:1.1}). We can estimate the
values of the parameters $p$ and $q$ there from given data of $V(\cdot )$ 
by
a least-squares approach (\cite{AIP}). In fact, since $V(\cdot )$ is a
stationary increment process, the variance of $V(t)-V(s)$ is a function in 
$%
t-s$. To be precise, 
\begin{equation*}
\frac{1}{t-s}\mathrm{Var}(V(t)-V(s))=U(t-s)\qquad (0\leq s<t),
\end{equation*}
where the function $U(t)=U(t;p,q)$ is given by 
\begin{equation}
U(t)=\frac{q^{2}}{(p+q)^{2}}+\frac{p(2q+p)}{(p+q)^{3}}\cdot \frac{%
(1-e^{-(p+q)t})}{t}\qquad (t>0).  \label{eq:5.1}
\end{equation}
Suppose that for $t_{j}=j\Delta t$, $j=1,\dots ,N$, the value of $V(t_{j})$
is $v_{j}$. For simplicity, we assume that $\Delta t=1$. The unbiased
estimation $u_{j}$ that corresponds to $U(t_{j})$ is given by 
\begin{equation*}
u_{j}=\frac{1}{j(N-j-1)}\sum_{i=1}^{N-j}(v_{i+j}-v_{i}-m_{j})^{2},
\end{equation*}
where $m_{j}$ is the mean of $v_{i+j}-v_{i}$'s: 
\begin{equation*}
m_{j}=\frac{1}{N-j}\sum_{i=1}^{N-j}(v_{i+j}-v_{i}).
\end{equation*}
Fitting $\{U(t_{j};p,q)\}$ to $\{u_{j}\}$ by least squares, we obtain the
desired estimated values of $p$ and $q$.

\begin{figure}[hbtp]
\begin{center}
\includegraphics[width=350pt,height=260pt]{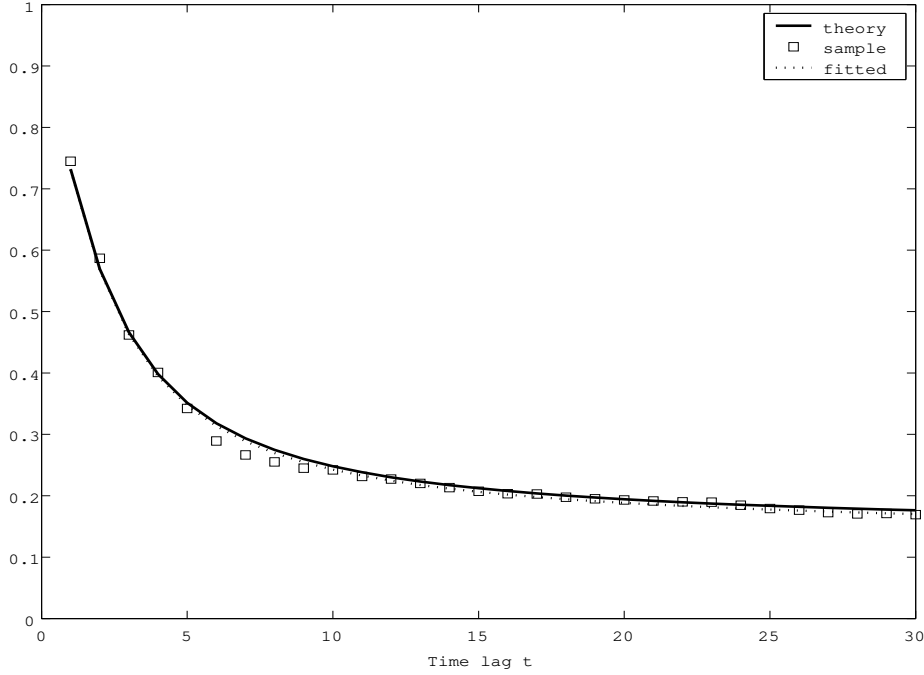}
\end{center}
\caption{Plotting of the function $v(t)$, the sample data and the fitted
function $v_0(t)$.}
\label{fig:5.1}
\end{figure}

For example, we produce sample values $v_{1},v_{2},\dots ,v_{1000}$ for $%
V(\cdot )$ with $(p,q)=(0.5,0.3)$ by a Monte Carlo simulation. We use this
data to numerically calculate the values $p_{0}$ and $q_{0}$ of $p$ and 
$q$,
respectively, that minimize 
\begin{equation*}
\sum_{j=1}^{30}(U(t_{j};p,q)-u_{j})^{2}.
\end{equation*}
It turns out that $p_{0}=0.5049$ and $q_{0}=0.2915$. In Figure 
\ref{fig:5.1}%
, we plot $\{U(t_{j};p,q)\}$ (theory), $\{u_{j}\}$ (sample) and $%
\{U(t_{j};p_{0},q_{0})\}$ (fitted). It is seen that the fitted curve 
follows
the theoretical curve reasonably well.

We extend the approach above to that for the estimation of the parameters 
$p$%
, $q$, $\theta$ and $\sigma$ in 
\begin{equation*}
dX(t)= -\theta X(t)dt + \sigma dV(t), \quad X(0)=X_0,
\end{equation*}
where $\theta,\sigma\in (0,\infty)$, the process $V(\cdot)$ is given by 
(\ref
{eq:1.1}) as above, and the initial value $X_0$ is independent of $%
(V(t))_{0\le t\le T}$ and satisfies $\mathbb{E}[X_0^2]<\infty$. The 
solution 
$X(\cdot)$ is the following Ornstein-Uhlenbeck-type process with memory: 
\begin{equation}
X(t)= e^{-\theta t}X_0 + \int_0^t e^{-\theta (t-u)}dV(u)\qquad (t\in 
[0,T]).
\label{eq:5.2}
\end{equation}
Put 
\begin{equation*}
\varphi(t):= \int_0^t e^{(\theta-p-q)u}du \qquad (0\le t\le T).
\end{equation*}

\begin{prop}
\label{prop:5.1} We have 
\begin{equation}
\frac{1}{t-s}\mathbb{E}\left[(X(t)-e^{-\theta(t-s)}X(s))^2\right]
=H(t-s)\quad (0\le s< t\le T),  \label{eq:5.3}
\end{equation}
where, for $0< t\le T$, the function $H(t)=H(t;p,q,\theta,\sigma)$ is given
by 
\begin{equation*}
H(t)=\sigma^2\left\{1-\frac{p(2q+p)}{(p+q)(\theta+p+q)}\right\} \frac{%
1-e^{-2\theta t}}{2\theta t} +\frac{\sigma^2 p(2q+p)e^{-2\theta 
t}\varphi(t)%
} {(p+q)(\theta+p+q)t}.
\end{equation*}
\end{prop}

\begin{proof}
By (\ref{eq:5.2}), the left-hand side of (\ref{eq:5.3}) is equal to
\begin{equation*}
\frac{\sigma^2 e^{-2\theta t}}{t-s}
  \mathbb{E}\left[\left(\int_s^t e^{\theta u}dV(u)\right)^2\right].
\end{equation*}
We put $c_u=pe^{-(p+q)u}I_{(0,\infty)}(u)$ for $u\in\mathbf{R}$. 
By Proposition 3.2 in \cite{AIK}, $\int_s^t e^{\theta u}dV(u)$ is given by
\begin{align*}
 &\int_s^t e^{\theta u}dW(u)
 -\int_{-\infty}^{\infty} \left(\int_s^t e^{\theta r}c_{r-u}
 dr\right)dW(u) \\
 &= \int_s^t\left(e^{\theta u}-\int_u^t e^{\theta r}c_{r-u}dr\right)dW(u)
 -\int_{-\infty}^s\left(\int_s^t e^{\theta r}c_{r-u}dr
  \right)dW(u).  
\end{align*} 
Thus $\mathbb{E}[(\int_s^t e^{\theta u}dV(u))^2]$ is equal to
\begin{equation*}
 \int_s^t \left(e^{2\theta u}-
   \int_u^t e^{\theta r}c(r-u)dr\right)^2du 
 +\int_{-\infty}^s \left(\int_s^t e^{\theta r}c(r-u)dr\right)^2du.  
\end{equation*}
By integration by parts and the equalities  
\begin{gather*}
 \int_u^t e^{\theta r}c_{r-u}dr= pe^{(p+q)u}\{\varphi(t)-\varphi(u)\},\\ 
 e^{-(\theta-p-q)s}\{\varphi(t)-\varphi(s)\}=\varphi(t-s),  
\end{gather*}
we obtain the desired result.
\end{proof}

Suppose that for $t_{j}=j\Delta t$, $j=1,\dots ,N$, the value of $X(t_{j})$
is $x_{j}$. We assume $\Delta t=1$ for simplicity. The estimation $%
h_{j}(\theta )$ that corresponds to $H(t_{j};p,q,\theta ,\sigma )$ is given
by 
\begin{equation*}
h_{j}(\theta )=\frac{1}{j(N-j-1)}\sum_{i=1}^{N-j}(x_{i+j}-e^{-\theta
j}x_{i}-m_{j}(\theta ))^{2},
\end{equation*}
where $m_{j}(\theta )$ is the mean of $x_{i+j}-e^{-\theta j}x_{i}$, $%
i=1,\dots ,N-j$: 
\begin{equation*}
m_{j}(\theta )=\frac{1}{N-j}\sum_{i=1}^{N-j}(x_{i+j}-e^{-\theta j}x_{i}).
\end{equation*}
Fitting $\{H(t_{j};p,q,\theta ,\sigma )-h_{j}(\theta )\}$ to $\{0\}$ by
least squares, we obtain the desired estimated values of $p$, $q$, $\theta 
$
and $\sigma $.

\begin{figure}[hbtp]
\begin{center}
\includegraphics[width=350pt,height=260pt]{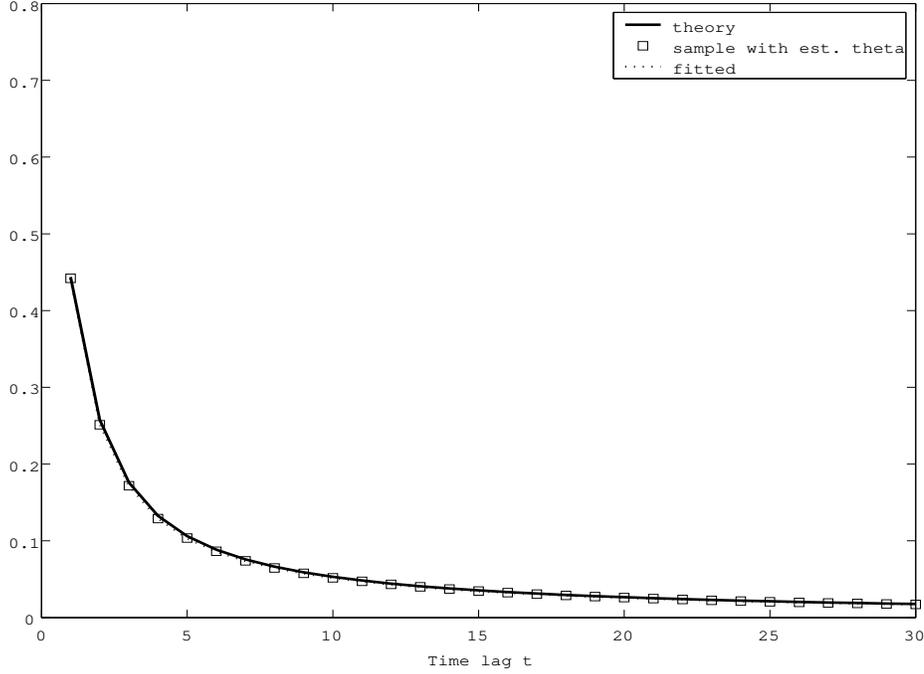}
\end{center}
\caption{Plotting of the function $h(t)$, the sample data $\bar{h}(t)$ with
estimated $\protect\theta$ and the fitted function $h_0(t)$.}
\label{fig:5.2}
\end{figure}

For example, we produce sample values $x_{1},x_{2},\dots ,x_{1000}$ for $%
X(\cdot )$ with 
\begin{equation*}
(p,q,\theta ,\sigma )=(0.2,1.5,0.8,1.0)
\end{equation*}
by a Monte Carlo simulation. We use this data to numerically calculate the
values $p_{0}$, $q_{0}$, $\theta _{0}$ and $\sigma _{0}$ of $p$, $q$, $%
\theta $ and $\sigma $, repectively, that minimize 
\begin{equation*}
\sum_{j=1}^{30}(H(t_{j};p,q,\theta ,\sigma )-h_{j}(\theta ))^{2}.
\end{equation*}
It turns out that 
\begin{equation*}
(p_{0},q_{0},\theta _{0},\sigma _{0})=(0.1910,1.5382,0.8354,1.0184).
\end{equation*}
In Figure \ref{fig:5.2}, we plot $\{H(t_{j};p,q,\theta ,\sigma )\}$
(theory), $\{h_{j}(\theta _{0})\}$ (sample with estimated~$\theta $) and $%
\{H(t_{j};p_{0},q_{0},\theta _{0},\sigma _{0})\}$ (fitted). It is seen that
the fitted curve follows closely the theoretical curve.

\section{Simulation}

\label{sec:6}

As we have seen, the process $V(\cdot )$ described by (\ref{eq:1.1}) has
both stationary increments and a simple semimartingale representation as
Brownian motion does, and it reduces to Brownian motion when $p=0$. In this
sense, we may see $V(\cdot )$ as a generalized Brownian motion. Since $%
V(\cdot )$ is non-Markovian unless $p=0$, we have now a wide choice for
designing models driven by either white or colored noise.

In this section, we compare the performance of the optimal filter with the
Kalman-Bucy filter in the presence of colored noise. We consider the
partially observable process $(X(t),Y(t))_{0\leq t\leq T}$ governed by 
(\ref
{eq:1.5}) with $X_{0}=0$. Suppose that an engineer uses the conventional
Markovian model 
\begin{equation*}
dX^{\prime }(t)=\theta X^{\prime }(t)dt+\sigma dB_{1}(t),\quad X^{\prime
}(0)=0
\end{equation*}
to describe the non-Markovian system process $X(\cdot )$. Then he will be
led to use the Kalman-Bucy filter $\widetilde{X}(\cdot )$ governed by 
\begin{equation*}
d\widetilde{X}(t)=(\theta -\mu ^{2}\gamma (t))\widetilde{X}(t)dt+\mu \gamma
(t)dY(t),\quad \widetilde{X}(0)=0,
\end{equation*}
where $\gamma (\cdot )$ is the solution to 
\begin{equation*}
\frac{d\gamma (t)}{dt}=\sigma ^{2}+2\theta \gamma (t)-\mu ^{2}\gamma
(t)^{2},\quad \gamma (0)=0
\end{equation*}
(see (\ref{eq:3.8}) and (\ref{eq:3.9})), instead of the right optimal 
filter 
$\widehat{X}$ as described by Theorem \ref{thm:3.4}.

We adopt the following parameters: 
\begin{gather*}
T=10,\quad \Delta t=0.01,\quad N=T/\Delta t=1000, \\
t_{i}=i\Delta t\quad (i=1,\dots ,N), \\
\sigma =1,\quad \theta =-2,\quad \mu =5.
\end{gather*}
Let $n\in \{1,2,\dots ,100\}$. For the $n$-th run of Monte Carlo 
simulation,
we sample $x_{n}(t_{1}),\dots ,x_{n}(t_{N})$ for $X(\cdot )$. Let $%
\widetilde{x}_{n}(\cdot )$ and $\widehat{x}_{n}(\cdot )$, $n=1,\dots ,100$,
be the Kalman-Bucy filter and the optimal filter respectively. For $u_{n}=%
\widehat{x}_{n}$ or $\widetilde{x}_{n}$, we consider the \textit{average
error norm\/} 
\begin{equation*}
\mathrm{AEN}:=\sqrt{\frac{1}{100N}\sum_{i=1}^{N}%
\sum_{n=1}^{100}(x_{n}(t_{i})-u_{n}(t_{i}))^{2}},
\end{equation*}
and the \textit{average error\/} 
\begin{equation*}
\mathrm{AE}(t_{i}):=\sqrt{\frac{1}{100}%
\sum_{n=1}^{100}(x_{n}(t_{i})-u_{n}(t_{i}))^{2}}\qquad (i=1,\dots ,N).
\end{equation*}

\begin{figure}[hbtp]
\begin{center}
\includegraphics[width=350pt,height=260pt]{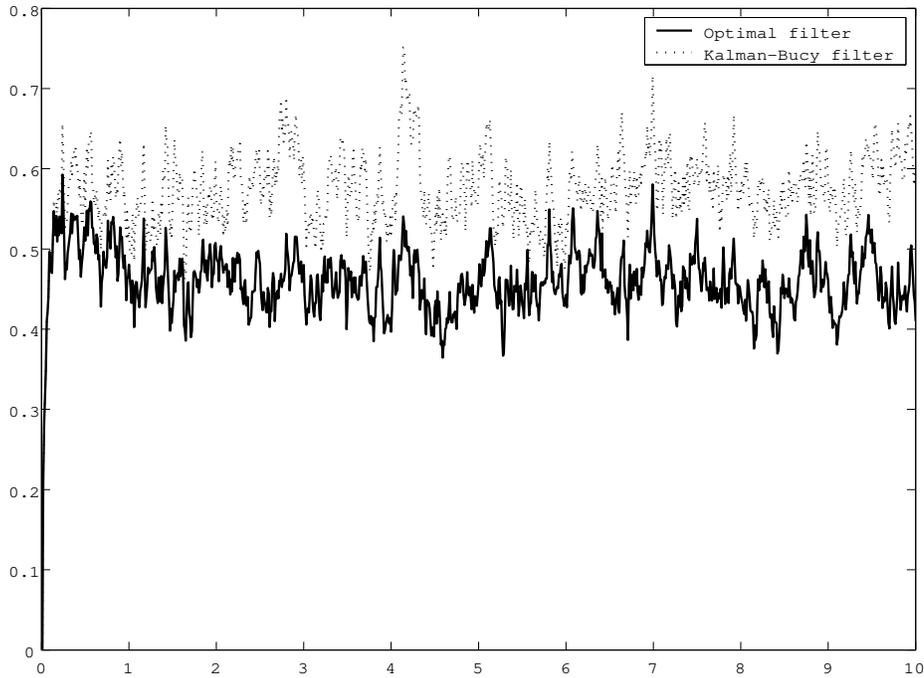}
\end{center}
\caption{Plotting of $\mathrm{AE}(\cdot)$ for the optimal and Kalman-Bucy
filters with noise parameter $\Theta=\Theta_2$}
\label{fig:6.1}
\end{figure}

In Table \ref{table:6.1}, we show the resulting $\mathrm{AEN}$'s of $\{%
\widehat{x}_{n}\}$ and $\{\widetilde{x}_{n}\}$ for the following five sets
of $\Theta =(p_{1},q_{1},p_{2},q_{2})$: 
\begin{align*}
\Theta _{1}& =(0.2,0.3,0.5,0.2), \\
\Theta _{2}& =(5.2,0.3,-0.5,0.6), \\
\Theta _{3}& =(0.0,1.0,5.8,0.7), \\
\Theta _{4}& =(5.4,0.8,0.0,1.0), \\
\Theta _{5}& =(5.1,2.3,4.9,1.3).
\end{align*}
We see that there are clear differences between the two filters in the 
cases 
$\Theta _{2}$ and $\Theta _{4}$. We notice that, in these two cases, 
$p_{1}$
is large than the parameters $p_{2}$ and $q_{2}$. In Figure \ref{fig:6.1},
we compare the graphs of $\mathrm{AE}(\cdot )$ for the two filters in the
case $\Theta =\Theta _{2}$. It is seen that the error of the optimal filter
is consistently smaller than that of the Kalman-Bucy filter.

\begin{table}[tbph]
\begin{center}
\begin{tabular}{c|c|c}
$\Theta$ & Optimal 
filter & Kalman-Bucy filter \\ \hline
$\Theta_1$ & 0.5663 & 0.5667 \\ \hline
$\Theta_2$ & 0.4620 & 0.5756 \\ \hline
$\Theta_3$ & 0.5136 & 0.5167 \\ \hline
$\Theta_4$ & 0.4487 & 0.5196 \\ \hline
$\Theta_5$ & 0.4294 & 0.4524 \\ \hline
\end{tabular}
\end{center}
\caption{Comparison of AEN's}
\label{table:6.1}
\end{table}

\end{document}